\documentclass[12pt]{article}%23.06.10 complete version
\usepackage{amsmath}
\usepackage{amssymb}
\usepackage{cite}
\newsymbol \blackbox 1004
\newcommand{\eh}{\hfill}\newlength{\sperr}

\newenvironment{proof}{{\settowidth{\sperr}{\bf\rm
Proof}%
\par\addvspace{0.3cm}\noindent\parbox[t]{1.3\sperr}
{\bf\rm P\eh r\eh o\eh o\eh f\eh }%
}}{\nopagebreak\mbox{}
$\blackbox$\par\addvspace{0.3cm}}

\def\nn{\nonumber}

\def\g{\gamma}
\def\G{\Gamma}

\def\s{\sigma}
\def\la{\lambda}
\def\om{\omega}

\def\vp{\varphi}
\def\ve{\varepsilon}
\def\wh{\widehat}
\def\wt{\widetilde}
\def\ov{\overline}
\def\ora{\overrightarrow}
\def\p{\partial}

\def\BC{{\mathbb C}}
\def\BR{{\mathbb R}}
\def\BN{{\mathbb N}}
\def\clp{{\mathcal P}}
\def\cla{{\mathcal A}}

\def\cld{{\mathcal D}}
\def\mclo{{\mathcal O}}

\newtheorem{Pa}{Paper}[section]
\newtheorem{Tm}[Pa]{{\bf Theorem}}

\newtheorem{Cy}[Pa]{{\bf Corollary}}
\newtheorem{Rk}[Pa]{{\bf Remark}}

\newtheorem{Pn}[Pa]{{\bf Proposition}}

\title{On the factorisation formula for fundamental solutions in the inverse spectral transform}

\author{Alexander Sakhnovich}

\date{}
\begin{document}
\maketitle

\thanks{Fakult\"at f\"ur Mathematik, Universit\"at Wien,
\\
Nordbergstrasse 15, A-1090 Vienna,  Austria. \\
al$_-$sakhnov@yahoo.com}

\begin{abstract} 
A factorization formula for wave functions, which is basic in the inverse spectral transform approach 
to initial-boundary value problems, is proved in greater generality than before.
Applications follow. Related compatibility  questions for the GBDT version of   B\"acklund-Darboux 
transformation are treated too.
\end{abstract}



\section{Introduction} \label{intro}
\setcounter{equation}{0}
Zero curvature representation of the integrable nonlinear equations
is a well known approach (see \cite{AKNS, TF, Nov, ZM} and references in \cite{TF}), 
which was developed
soon after the seminal Lax pairs appeared in \cite{Lax}. Namely, many
integrable nonlinear equations admit representation
(zero curvature representation)
\begin{align}&      \label{0.1}
G_t(x,t,z)-F_x(x,t,z)+[G(x,t,z),F(x,t,z)]=0, \\
& \nonumber
 G_t:=\frac{\p}{\p t}G, \quad [G,F]:=GF-FG,
\end{align}
which is the compatibility condition of the auxiliary linear systems
\begin{equation}      \label{0.2}
\frac{\p}{\p x}w(x,t,z)=G(x,t,z)w(x,t,z), \quad \frac{\p}{\p t}w(x,t,z)=F(x,t,z)w(x,t,z).
\end{equation}
Here $G$ and $F$ are $m \times m$ matrix functions, and
$z$ is the spectral parameter, which will be omitted sometimes in our notations.

Solution of integrable nonlinear equations is closely related to Lax pairs
and  zero curvature representations, which have been mentioned above, and has been
a great breakthrough in the second half of the 20th century. An active study
of the cases, which are close to integrable in a certain sense,  followed
(see, for instance, some references in \cite{Bo, Persp, Ka}).
Initial-boundary value problems  for  integrable nonlinear equations
can be considered as an important example, where integrability is
"spoiled" by the boundary conditions. These problems are
of great current interest, and inverse spectral transform (ISpT) method
\cite{Ber, BerG, KvM, SaA1, SaA2, SaA4, SaA5, SaL1, SaL2, SaL3, SaL4}
is one of the fruitful approaches in this domain. Further we assume that
$x, \, t$ belong to a semi-strip
\begin{equation}      \label{0.3}
{\mathcal D}=\{(x,\, t):\,0 \leq x <\infty, \,\,0\leq t<a\}.
\end{equation}
Normalize fundamental solutions of the auxiliary systems by the initial conditions
\begin{align} \label{0.4}
&\frac{d}{d x}W(x,t,z)=G(x,t,z)W(x,t,z), \quad W(0,t,z)=I_m;
\\ \label{0.5} &
\frac{d}{d t}R(x,t,z)=F(x,t,z)R(x,t,z), \quad R(x,0,z)=I_m,
\end{align}
where $I_m$ is the identity matrix of order $m$.
If condition \eqref{0.1} holds, the fundamental solution  of \eqref{0.4}
 admits factorization
\begin{equation} \label{0.6}
W(x,t,z)=R(x,t,z)W(x,0,z)R(t,z)^{-1}, \quad R(t,z):=R(0,t,z).
\end{equation}
Formula \eqref{0.6} is one of the basic and actively used formulas in the
 inverse spectral transform method (see \cite{SaA1, SaA2, SaA4, SaA5, SaL1, SaL2, SaL3, SaL4}
and references therein). It was derived in \cite{SaL1, SaL2} under some smoothness conditions
(conitinuous differentiability of $G$ and $F$, in particular): see formulas (1.6) in \cite{SaL1}, p.22
and in \cite{SaL2}, p. 39.

Here we prove \eqref{0.6}  under weaker conditions and in much greater detail, which is important
for applications.  Namely, we prove the following theorem.
\begin{Tm} \label{TmM} Let $m \times m$ matrix functions $G$ and $F$ and their derivatives $G_t$ and
$F_x$  exist on the semi-strip $\mathcal{D}$,
let $G$, $G_t$, and $F$ be continuous with respect to $x$ and $t$ on $\mathcal{D}$,
 and let \eqref{0.1}
hold. Then the equality 
\begin{equation} \label{0.6'}
W(x,t,z)R(t,z)=R(x,t,z)W(x,0,z), \quad R(t,z):=R(0,t,z),
\end{equation}
is true.
 \end{Tm}
 Note that constructions similar to \eqref{0.6} appear also in the theory of
 Knizhnik-Zamolodchikov equation (see Theorem 3.1 in \cite{SaLKZ2} and see also \cite{SaLKZ}).

Theorem \ref{TmM} is proved in Section \ref{Pr}.
Section \ref{some} is dedicated to  applications to initial-boundary value
problems, and Theorem \ref{evol} on the evolution of the Weyl function
for the "focusing"  modified Korteweg-de Vries (mKdV) equation is proved
there as an example.

Related questions of the equality of mixed derivatives and application
of this equality  to the GBDT version (see \cite{GKS2, GKS6, KaSa, MST, SaA2, SaA2', SaAJFA, SaA3, SaA6}
and references therein) 
of the B\"acklund-Darboux transformation are 
treated in Section \ref{Darb}.

As usual, by $\BN$ we denote the set of positive integers, by $\BC$ we denote the complex plane, 
and by $\BC^m$ is denoted the $m$-dimensional coordinate space 
over $\BC$.  By $\Im z$ is denoted  the imaginary part of $z \in \BC$,
and $\arg z$ is the argument of $z$.
By $C^k({\mathcal D})$ we denote functions and matrix functions, which are $k$ times continuously differentiable on ${\cal D}$.

\section{Proof of Theorem \ref{TmM}} \label{Pr}
\setcounter{equation}{0}
The spectral parameter $z$ is non-essential for the formulation of 
Theorem \ref{TmM} and for its proof and we shall omit it in this section.
We shall need the proposition below.
\begin{Pn} \label{PM} Let the $m \times m$ matrix function $W$ be given on the
 semi-strip $\mathcal{D}$ by  equation \eqref{0.4}, where $G(x,t)$
and $G_t(x,t)$ are continuous matrix functions in $x$ and $t$.

(i)Then the derivative $W_t$ exists and matrix functions $W$ and $W_t$
are continuous with respect to $x$ and $t$
on the semi-strip $\mathcal{D}$.

(ii)Moreover, the mixed derivative $W_{tx}$ exists and the equality
$W_{tx}=W_{xt}$ holds on $\mathcal{D}$.
 \end{Pn}
\begin{proof}. 
Consider system
 \begin{align}      \label{1.2}
&\frac{d}{d x}y=\wh G(x,y)y, \quad
 \wh G(x,y)=\wh G(x,y_{m+1}):=\left[
\begin{array}{cc}
G(x,y_{m+1}) &0\\ 0&0
\end{array}
\right],  
\end{align}
where $\wh G$ is an $m+1 \times m+1$ matrix function and $y_{m+1}$ is the last entry of the column vector $y\in \BC^{m+1}$.
Denote by $W_j$ and $e_j$ the $j$-th columns of $W$ and $I_m$,
respectively ($1 \leq j \leq m$).
 It easily follows from \eqref{0.4} that the solution of \eqref{1.2} with the initial condition
 \begin{align}      \label{1.3}
& y(0)=g=\left[
\begin{array}{c}e_j \\ t
\end{array}
\right]
\end{align}
has the form
\begin{equation}      \label{1.4}
y(x,g)=\left[
\begin{array}{c}W_j(x,t) \\ t
\end{array}
\right].
\end{equation}
Putting $G(x,t)=G(0,t)$ for $-\ve \leq x \leq 0$ whereas $t \geq 0$, and putting $G(x,t)=G(x,0)+tG_t(x,0)$
for $-\ve \leq t \leq 0$ ($\ve>0$) we extend $G$ so that $G$ and $G_t$ remain continuous on the
rectangles
\begin{equation}      \label{1.5}
\mathcal{D}(a_1,\, a_2)=\{(x,\, t): \,  -\ve\leq x \leq a_1,  \, \, -\ve \leq t \leq a_2<a\},
\quad a_1, a_2 \, \in \BR_+.
\end{equation}
Hence, it follows from the definition
of $\wh G$ in \eqref{1.2} that $\wh G(x,y)$ and, as a consequence, the vector function $\wh G(x,y)y$ are continuous  on $\mathcal{D}(a_1,\, a_2)$ together
with their derivatives  with respect to the entries of $y$. Thus, according  to the classical theory
of ordinary differential equations (see, for instance, theorem on pp. 305-306 in \cite{Sm})
the
partial first derivatives of $y(x,g)$ with respect to the entries of $g$
exist in the interior
$\mathcal{D}_i(a_1,\, a_2)$ of  $\mathcal{D}(a_1,\, a_2)$.
Moreover,
$y$ and its partial  derivatives  with respect to the entries of $g$ are continuous.  In particular,  
since by
 \eqref{1.3} we have $g_{m+1}=t$,
the functions $y$ and $y_t$ are continuous in all rectangles $\mathcal{D}_i(a_1,\, a_2)$.
Taking into account \eqref{1.4}, we see that $W$ and $W_t$
are continuous  in the rectangles $\mathcal{D}_i(a_1,\, a_2)$, and the statement (i) is true.

In view of \eqref{0.4} and considerations above the derivatives $W_x$, $W_{xt}$, and $W_t$ exist and are
continuous in the rectangles $\mathcal{D}_i(a_1,\, a_2)$. Hence, by a  stronger formulation (see, for instance, 
\cite{See, Aks} or p. 201 in \cite{Mar}) of the well-known theorem on mixed derivatives, $W_{tx}$ exists in 
$\mathcal{D}_i(a_1,\, a_2)$  and $W_{tx}=W_{xt}$. Thus,  the statement (ii) follows.
\end{proof}
Now, we can follow the scheme from Chapter 3 in \cite{SaL2}  (see also Chapter  12 in \cite{SaL4}).
\begin{proof} of  Theorem \ref{TmM}.
According to statement (i) in Proposition  \ref{PM} the matrix function $W_t$ exists
and is continuous. Introduce $U(x,t)$ by the equality
\begin{equation}      \label{1.6}
U:=W_t-FW.
\end{equation}
By \eqref{0.4},  \eqref{1.6}, and  statement (ii) in Proposition  \ref{PM}  we have
\begin{equation}      \label{1.7}
U_x=W_{tx}-F_xW-FW_x=W_{xt}-F_xW-FGW.
\end{equation}
It  is immediate also from \eqref{0.4}  that
\begin{equation}      \label{1.8}
 W_{xt}=\big(GW\big)_t=G_tW+GW_t.
\end{equation}
Formulas  \eqref{1.7} and   \eqref{1.8} imply
\begin{equation}      \label{1.9}
 U_{x}=G_tW+GW_t-F_xW-FGW=(G_t-F_x+GF-FG)W+GW_t-GFW.
\end{equation}
It follows from \eqref{0.1}, \eqref{1.9}, and definition  \eqref{1.6}  that $U_{x}=GU$,
that is, $U$ and $W$ satisfy the same equation. Taking into account $W(0,t)=I_2$,
we derive $W_t(0,t)=0$, and so by \eqref{1.6} we have $U(0,t)=-F(0,t)$.  Finally, as 
\[
U_{x}=GU, \quad W_x=GW, \quad U(0,t)=-F(0,t), \quad W(0,t)=I_2,
\]
we have $U(x,t)=-W(x,t)F(0,t)$ or, equivalently,
\begin{equation}      \label{1.10}
W_t(x,t)-F(x,t)W(x,t)= -W(x,t)F(0,t).
 \end{equation}
 Put
 \begin{equation}      \label{1.11}
Y(x,t)=W(x,t)R(t), \quad Z(x,t)=R(x,t)W(x,0).
 \end{equation}
 Recall that $R(t)=R(0,t)$. Therefore 
 \eqref{0.5}, \eqref{1.10}, and \eqref{1.11}  imply that
\begin{align}     \nonumber
Y_t(x,t)&=\big(F(x,t)W(x,t)-W(x,t)F(0,t)\big)R(t)+W(x,t)F(0,t)\big)R(t)
\\ &  \label{1.12}
=F(x,t)Y(x,t), \quad Y(x,0)=W(x,0).
 \end{align}
 Formulas \eqref{0.5} and \eqref{1.11}  imply that
 \begin{equation}      \label{1.13}
Z_t(x,t)=F(x,t)Z(x,t), \quad Z(x,0)=W(x,0).
 \end{equation}
 By \eqref{1.12} and \eqref{1.13} $Y=Z$, that is, \eqref{0.6'} holds.
\end{proof}
\begin{Rk} \label{Dom0}
Though the case of continuous $F$ is more convenient for applications, 
it is immediate from the proof that the statement of Theorem  \ref{TmM}
is true, when $F$ is differentiable with respect to $x$, and measurable
and summable  with respect to $t$ on all  finite intervals from $\BR_+$ .
\end{Rk}
According to the proof of 
Theorem \ref{TmM} the following remark is also true.
\begin{Rk} \label{Dom}
Theorem \ref{TmM} holds on the domains more general than $\cal D$.
In particular, it holds
if we consider  $(x,\, t)\in {\cal I}_1 \times {\cal I}_2$,
where ${\cal I}_k$
 $(k=1,2)$  is the interval $[0, \,b_k)$  $\,(0<b_k\leq \infty)$.
\end{Rk}

Another interesting case of matrix factorizations related to boundary value
problems is treated in \cite{BHLW, KL}.

%%%%%%%%%%%%%%%%%%%%%%%%%%%%%%%%%%%%%%%%%%%%
%%%%%%%%%%%%%%%%%%%%%%%%%%%%%%%%%%%%%%%%%%%%
\section{Some applications} \label{some}
\setcounter{equation}{0}
The {\it matrix "focusing" mKdV equation} has the form
\begin{align}  \label{2.1} &
4v_t=v_{xxx}+3\big(v_xv^*v+vv^*v_x\big),
\end{align}
where $v(x,t)$ is a $p \times p$ matrix function.
 Equation \eqref{2.1} is equivalent 
(see \cite{CD, TF, Wad} and references therein)
to zero curvature
equation \eqref{0.1}, where the $m \times m$ ($m=2p$)
matrix functions $G(x,t,z)$ and $F(x,t,z)$ are given by the formulas
\begin{align}  \label{2.2} &
G=izj+V, \quad
j=\left[
\begin{array}{lr}
I_p& 0 \\ 0 & -I_p
\end{array}
\right],
\quad
 V=\left[
\begin{array}{lr}
0& v \\ - v^* & 0
\end{array}
\right],
\\ \label{2.3} &
F=-iz^3j-z^2V-\frac{iz}{2}\big(V^2+V_xj\big)+\frac{1}{4}
\big(V_{xx}-2V^3-V_xV+VV_x\big).
\end{align}
At first we omit the variable $t$ in $V$ and $v$.
The Weyl theory of  the skew-self-adjoint Dirac system (also called
Zakharov-Shabat or AKNS system)
 \begin{equation}      \label{2.4}
\frac{d}{dx}w(x,z)=\big(izj+V(x)\big)w(x,z), \quad x \geq 0
 \end{equation}
 was treated in \cite{CG, FKS, SaA1, SaA2} (see also preliminaries
 in \cite{SaA7}). 
 
 For the case of measurable matrix function $v$ such that
 \begin{equation} \label{2.5}
\sup_{0<x<\infty}\|v(x)\| \leq M,
\end{equation}
the Weyl matrix function $\vp$ of system \eqref{2.4} is  uniquely  defined in the semi-plane
$\Im z<-M$ by the inequality
\begin{equation} \label{2.6}
\int_0^\infty \left[ \begin{array}{lr}   \varphi (z)^* &
I_p
\end{array} \right]
  W(x, z)^*
W(x, z)
\left[ \begin{array}{c}
  \varphi (z) \\ I_p \end{array} \right] dx < \infty, \quad  \Im z<-M<0,
\end{equation}
where $W$ is the normalized by $W(0,z)=I_m$  fundamental solution of \eqref{2.4}.
Weyl functions are constructed using pairs of
meromorphic $p \times p$  matrix functions $P_1(z), \, P_2(z)$,  which are nonsingular
and have property-$j$, that is,
 \begin{equation} \label{2.6'}
\clp(z)^*\clp(z)>0, \quad \clp(z)^*j\clp(z) \leq 0, \quad \clp: \left[
\begin{array}{c}
P_1 \\ P_2
\end{array}
\right].
\end{equation}
 \begin{Tm}\label{DP} \cite{SaA1} There is a unique Weyl function of such a system
 \eqref{2.4}  that \eqref{2.5} holds. This Weyl function is holomorphic
 in the semi-plane $\Im z<-M$.
 It is
 given by the equality
 \begin{align} \label{2.7}
\vp(z)=& \lim_{r \to \infty}\big({\cal A}_{11}(r,z)P_1(r,z)+{\cal A}_{12}(r,z)P_2(r,z)\big)
\\ \nonumber & \times
\big({\cal A}_{21}(r,z)P_1(r,z)+{\cal A}_{22}(r,z)P_2(r,z)\big)^{-1} \quad (\Im z<-M),
\\ \label{2.8} &
{\cal A}(r,z)=\big\{{\cal A}_{kp}(r,z)  \big\}_{k,p=1}^2:=W(r,\ov z)^*,
\end{align}
 where the pairs $\{P_1, \, P_2\}$ are arbitrary  pairs satisfying \eqref{2.6'}.
 \end{Tm}
 Our next theorem on the evolution of the Weyl function
 in the case of the focusing mKdV follows from Theorems \ref{TmM} and \ref{DP}. 
 The case of the defocusing mKdV was earlier treated in \cite{SaL1, SaL2, SaL4}.
  \begin{Tm}\label{evol} Let a $p \times p$ matrix function $v \in C^1({\mathcal D})$ 
have a continuous partial second derivative  $v_{xx}$, and let $v_{xxx}$  exist.
 Assume that $v$  satisfies mKdV \eqref{2.1}  and that the inequalities
  \begin{equation} \label{2.8'}
  \sup_{(x,t)\in {\mathcal D}}\|v(x,t)\| \leq M, \quad
\sup_{(x,t)\in {\mathcal D}}\big( \|v_x(x,t)\|+ \|v_{xx}(x,t)\| \big) <\infty 
\end{equation}
hold.

 Then the evolution $\vp(t,z)$ of the Weyl function  of the skew-self-adjoint Dirac system
 \eqref{0.4}, where $G$ has the form \eqref{2.2}, is given by the equality
 \begin{equation} \label{2.9}
\vp(t,z)=\big(R_{11}(t,z)\vp(0,z)+R_{12}(t,z)\big)
\big(R_{21}(t,z)\vp(0,z)+R_{22}(t,z)\big)^{-1}
\end{equation}
in the semi-plane $\Im z <-M<0$.
Here the block matrix function 
 \begin{equation} \label{2.9'}
R(t,z)=\{R_{kn}(t,z)\}_{k,n=1}^2=R(0,t,z)
\end{equation}
is defined
by the boundary values $v(0,t)$, $v_x(0,t)$, and $v_{xx}(0,t)$ via formulas
\eqref{0.5} and \eqref{2.3}.
  \end{Tm}
  \begin{proof}. As $V^*=-V$ and $(V_x j)^*=V_xj$, it is immediate from \eqref{2.3}
  that $F(x,t, \ov z)^*+F(x,t,z)=0$. Hence, it follows from
  \eqref{0.5}  that 
\[  
  \frac{\p}{\p t}\big(R(x,t, \ov z)^*R(x,t,z)\big)=0.
\]  
  Therefore, using equalities  $R(x,0,z)=I_m$ and \eqref{2.9'}, we get
 \[
 R(x,t, \ov z)^*R(x,t,z)=I_m, \qquad R(t, \ov z)^*R(t,z)=I_m,
 \]
 or, equivalently,
\begin{equation} \label{2.10}
R(x,t, \ov z)^*=R(x,t,z)^{-1}, \qquad R(t, \ov z)^*=R(t,z)^{-1}.
\end{equation}
In view of  \eqref{2.10} rewrite \eqref{0.6'} in the form
\begin{equation} \label{2.11}
{\cal A}(x,t, z)R(x,t,z)=R(t,z){\cal A}(x,0, z),
\end{equation}
where  ${\cal A}(x,t, z):=W(x,t, \ov z)^*$ (compare with \eqref{2.8}).
Let ${\cal P}(x,z)$ satisfy \eqref{2.6'} and put
\begin{equation} \label{2.12}
\wt\clp(x,t,z)=\left[
\begin{array}{c}
\wt P_1(x,t,z) \\ \wt P_2(x,t,z)
\end{array}
\right]=R(x,t,z)\clp(x,z).
\end{equation}
By \eqref{2.11} and \eqref{2.12} we have
\begin{equation} \label{2.12'}
{\cal A}(x,t, z)\wt\clp(x,t,z)=R(t,z){\cal A}(x,0, z)\clp(x,z).
\end{equation}

Now, taking into account that $\clp(x,z)$ is a nonsingular pair with property-$j$,
we show that $\wt\clp(x,t,z)$ is a nonsingular pair with property-$j$ too. 
According to  \eqref{0.5}, \eqref{2.3}, and \eqref{2.8'} we get
\begin{equation} \label{2.13}
 \frac{\p}{\p t}\big(R(x,t, z)^*jR(x,t,z)\big)=R(x,t, z)^*\big(i(\ov z^3-z^3)I_m+O(z^2)\big)R(x,t,z)
\end{equation}
for $z \to \infty$. Formula \eqref{2.13} implies that for some 
\begin{equation} \label{2.14}
M_1>M>0 \quad (M \geq \sup_{(x,t)\in {\mathcal D}}\|v(x,t)\|),
\end{equation}
and for all $z$ from the domain
\begin{equation} \label{2.15}
{\mathcal D}_1=\{z: \,z\in \BC, \quad \Im z <-M_1,  \quad 0>\arg z > - \pi/4\}
\end{equation}
 we have
 \[
  \frac{\p}{\p t}\big(R(x,t, z)^*jR(x,t,z)\big) \leq 0,
 \]
 and so
 \begin{equation} \label{2.16}
 R(x,t, z)^*jR(x,t,z)\leq j.
 \end{equation}
 Relations \eqref{2.6'}, \eqref{2.12}, and \eqref{2.16} imply that
  \begin{equation} \label{2.17}
\wt \clp(x,t,z)^* \wt \clp(x,t,z)>0, \quad \wt \clp(x,t,z)^*j \wt \clp(x,t,z) \leq 0 \quad (z\in {\cal D}_1). 
\end{equation}
Clearly, it suffices to prove \eqref{2.9} for values of $z$ from ${\cal D}_1$. 
(According to \eqref{2.14} and \eqref{2.15} the domain ${\cal D}_1$ belongs to the semi-plane
$\Im z<-M$.)

In a way similar to the proofs of \eqref{2.10} and \eqref{2.16} we derive
 \begin{equation} \label{2.17'}
 \cla(x,t, z)=W(x,t,z)^{-1}, \qquad    W(x,t, z)^*jW(x,t,z)\geq j \quad (\Im z<-M).
 \end{equation}
It is immediate from \eqref{2.17'} that  
\begin{align} \label{o0}&
\cla(x,t, z)^*j\cla(x,t,z)\leq j \qquad (\Im z<-M).
 \end{align}
Hence, inequalities
\eqref{2.6'} and \eqref{2.17} imply
\begin{align} \label{o1}&
\det\big({\cal A}_{21}(x,0,z) P_1(x,z)+{\cal A}_{22}(x,0,z)P_2(x,z)\big)\not=0 \quad (\Im z<-M), \\
\label{o2}&
\det\big({\cal A}_{21}(x,t,z)\wt P_1(x,t,z)+{\cal A}_{22}(x,t,z)\wt P_2(x,t,z)\big)\not=0 
\quad (z\in {\cal D}_1).
 \end{align}
 In view of \eqref{o1} rewrite \eqref{2.12'} as
 \begin{align}\nonumber
{\cal A}(x,t, z)\wt\clp(x,t,z)=&R(t,z)
\left[
\begin{array}{c}
\phi(x,0,z) \\ I_p
\end{array}
\right]
\\ \label{o3} &
\big({\cal A}_{21}(x,0,z) P_1(x,z)+{\cal A}_{22}(x,0,z)P_2(x,z)\big), 
\end{align}
\begin{align}
\nonumber
\phi(x,0,z):=& \big({\cal A}_{11}(x,0,z) P_1(x,z)+{\cal A}_{12}(x,0,z)P_2(x,z)\big)
\\ & \label{o4} \times
\big({\cal A}_{21}(x,0,z) P_1(x,z)+{\cal A}_{22}(x,0,z)P_2(x,z)\big)^{-1}.
\end{align}
According to \eqref{o1}-\eqref{o3} we get
\begin{align} \nonumber&
\big({\cal A}_{11}(x,t,z) \wt P_1(x,t,z)+{\cal A}_{12}(x,t,z)\wt P_2(x,t,z)\big)
\\ & \label{o5}  \times
\big({\cal A}_{21}(x,t,z) \wt P_1(x,t, z)
+{\cal A}_{22}(x,t,z)\wt P_2(x,t,z)\big)^{-1}
\\ & \nonumber
=
\big(R_{11}(t,z) \phi(x,0,z)+R_{12}(t,z)\big)
\big(R_{21}(t,z) \phi(x,0,z)+R_{22}(t,z)\big)^{-1}.
\end{align}
As $\wt \clp(x,t,z)$ satisfies \eqref{2.17} for $z \in {\cal D}_1$, using \eqref{2.7} we derive
\begin{align} \nonumber
\vp(t,z)=& \lim_{x \to \infty} \big({\cal A}_{11}(x,t,z) \wt P_1(x,t,z)+{\cal A}_{12}(x,t,z)\wt P_2(x,t,z)\big)
\\ & \label{o6}  \times
\big({\cal A}_{21}(x,t,z) \wt P_1(x,t, z)
+{\cal A}_{22}(x,t,z)\wt P_2(x,t,z)\big)^{-1} \quad (z \in {\cal D}_1).
\end{align}
In a similar way we derive from \eqref{2.7} and  \eqref{o4} that
\begin{align}  & \label{o7}
\vp(0,z)=\lim_{x \to \infty}\phi(x,0,z) \quad (\Im z <-M).
\end{align}
Let us show that
 \begin{equation} \label{o8}
\det
\big(R_{21}(t,z)\vp(0,z)+R_{22}(t,z)\big)\not=0  \quad (z \in {\cal D}_1).
\end{equation}
Indeed, it follows from \eqref{2.6'}, \eqref{o0}, and  \eqref{o4} that
 \begin{align}\label{o9}
&
[\phi(x,0,z)^* \quad I_p] j
\left[
\begin{array}{c}
\phi(x,0,z) \\ I_p
\end{array}
\right]\leq 0.
\end{align}
By \eqref{o7} and \eqref{o9} the inequality
\begin{align}\label{o10}
&
[\vp(0,z)^* \quad I_p] j
\left[
\begin{array}{c}
\vp(0,z) \\ I_p
\end{array}
\right]\leq 0.
\end{align}
is true. Finally, inequalities \eqref{2.16} and \eqref{o10} imply
\begin{align}\label{o11}
&
[\vp(0,z)^* \quad I_p] R(t,z)^*jR(t,z)
\left[
\begin{array}{c}
\vp(0,z) \\ I_p
\end{array}
\right]\leq 0 \quad (z\in {\cal D}_1).
\end{align}
It is immediate from \eqref{o11} that  \eqref{o8} holds.
Relations \eqref{o5}-\eqref{o8} imply  \eqref{2.9} in the domain ${\cal D}_1$.
Hence, by analyticity  equality \eqref{2.9} holds  in the semi-plane $\Im z <-M$.
\end{proof}
In a way similar to \cite{Kr} and to more general constructions for self-adjoint systems in \cite{SaL2,
SaL4} (see also some references therein), one can use structured operators to solve inverse
problem for system \eqref{2.4} too. Namely, to recover $v$,  which satisfies condition \eqref{2.5},  
from the Weyl function $\vp$ we use operators $S_l$ (acting in $L_p^2(0,l)$, $0<l<\infty$) of the form
\begin{equation}\label{2.18}
S_lf=f(x)+ \frac{1}{2}\int_0^l
\int^{x+r}_{|x-r|} s^{\prime}\Big(
\frac{\la+x-r}{2} \Big) s^{\prime} \Big( \frac{\la+r-x}{2} \Big)^*d zf(r) dr.
\end{equation}
Here $s^{\prime}(x):
=\frac{d}{dx}s(x)$.  Below we give the procedure from \cite{SaA1} modified in accordance with \cite{FKS, SaA2}.

First, we recover
a $p \times p$ matrix function $s(x)$  with
the entries from $L^2(0,l)$ (i.e., $s(x) \in L^2_{p \times p}(0,l)$) via the
Fourier transform. That is, we put
\begin{equation}\label{2.19}
\displaystyle{s(x)= \frac{i}{2 \pi}e^{- \eta x}{\mathrm{
l.i.m.}}_{a \to \infty} \int_{- a}^{a}e^{i \xi
x} z^{-1} \vp(z /2) d \xi \quad (z= \xi +i \eta , \quad \eta
<-2M),}
\end{equation}
the limit l.i.m. being the limit in $L^2_{p \times p}(0,l)$.  Formula (\ref{2.19}) has sense for any $l<\infty$,
and so
the matrix function $s(x)$ is defined on the non-negative real  semi-axis $x\geq 0$.  Moreover,
 $s$ is absolutely continuous, it does not depend on the choice of $\eta <-
2M$, $s^{\prime}$ is bounded on any finite interval, and $s(0)=0$.
To define the operator $S_l$ we substitute $s^{\prime}(x)$
into  (\ref{2.18}).

Next, denote the $p\times 2p$ block rows of $W$ by $\om_1$ and $\om_2$:
\begin{equation}\label{2.20}
\om_1(x)=[ I_p \quad 0]W(x, 0), \quad  \om_2(x)=[0 \quad I_p]W(x, 0).
\end{equation}
It follows from (\ref{2.2}) and (\ref{2.4}) that $W(x,0)^*W(x,0)=I_{m}$. Hence, by 
(\ref{2.2}) , (\ref{2.4}),  and  (\ref{2.20}) we have
\begin{equation}\label{2.21}
v(x)=\om_1^{\prime}(x) \om_2(x)^*,
\end{equation}
and $\om_1$, $\om_2$ satisfy the equalities
\begin{equation} \label{2.22}
\om_1(0)=[I_p \quad 0], \quad \om_1 \om_1^* \equiv I_p, \quad
\om_1^{\prime}
\om_1^* \equiv 0, \quad 
\om_1 \om_2^* \equiv 0.
\end{equation}
It is immediate that $\om_1$ is uniquely recovered from $\om_2$ using  (\ref{2.22}).

Finally, we obtain $\om_2$ via the formula
\begin{equation}\label{2.23}
\om_2(l)=[0 \quad I_p]- \int_0^l \Big(S_l^{-1} s'(x)\Big)^*[I_p
\quad s(x)]dx \quad (0<l< \infty),
\end{equation}
where $S_l^{-1}$ is applied to $ s^{\prime}$ columnwise.
\begin{Tm}\label{IP}
Assume that $\vp$ is the Weyl function of system \eqref{2.4}, where $j$ and $V$ have the form
\eqref{2.2} and $v$ satisfies \eqref{2.5}. Then $v$ is recovered from $\vp$ via formulas
\eqref{2.21}-\eqref{2.23}, where $s$ and $S_l$ are given by equalities \eqref{2.18} and \eqref{2.19}.
All the mentioned above relations are well-defined and the inequalities $S_l \geq I$ hold.
\end{Tm}
Another inverse problem, where condition \eqref{2.5} on $v$ is substituted by a condition
on $\vp$, is also solved in \cite{FKS, SaA2, SaA7} using the same procedure.
\begin{Rk} \label{step} One can apply Theorems \ref{evol} and \ref{IP}
to recover solutions of mKdV. Theorems on the evolution of the Weyl functions
constitute also the first step in proofs of uniqueness and existence
of the solutions of nonlinear equations via  ISpT method (see, for instance, \cite{SaA7}).
\end{Rk}
%%%%%%%%%%%%%%%%%%%%%%%%%%%%%%%%%%%%%%%%%%%
\section{Factorization of the fundamental solution via Darboux matrix} \label{Darb}
\setcounter{equation}{0}
Various versions of B\"acklund-Darboux transformation and commutation \\
methods are widely used in spectral theory, differential equations
and nonlinear integrable equations (see, for instance, 
\cite{CC, Ci, Cr, D, Gu, GeT, Kr2, March, MS, ZM2}
and numerous references therein). In this section we consider a so called GBDT version
of the B\"acklund-Darboux transformation (see references in Introduction
and some basic notations and results in Appendix).
The statement of Theorem \ref{WF} is formulated and proved here in greater generality than before.

One can apply Theorem \ref{GBDTrd} on GBDT to construct solutions  and
  wave functions of nonlinear integrable equations.
For this purpose we use auxiliary linear systems for integrable nonlinear
equation, namely, linear systems :
\begin{eqnarray}       \label{g15}
&&w_x=Gw, \quad  w_t=Fw;
\\ && \label{g16}
G(x,t,z)=-\sum_{k=0}^r  z^k q_k(x,t)-\sum_{s=1}^l \sum_{k=1}^{r_s}(z-c_s)^{-k}
q_{sk}(x,t),
\\ && \label{g17}
F(x,t,z)=-\sum_{k=0}^R  z^k Q_k(x,t)-\sum_{s=1}^L \sum_{k=1}^{R_s}(z-C_s)^{-k}
Q_{sk}(x,t),
\end{eqnarray}
and zero curvature (compatibility condition) representation (\ref{0.1}) of the
integrable nonlinear equation itself. We consider nonlinear equations on the domain
$(x,t) \in  {\cal I}_1 \times {\cal I}_2$, where ${\cal I}_k$
 ($k=1,2$)   is the interval $[0, \,b_k)$,  $\,(0<b_k\leq \infty)$.
By  Theorem \ref{TmM} and Remark \ref{Dom} the following corollary is true.
\begin{Cy}\label{Rk3.3} 
Let coefficients $\{q_k(x,t)\}$ and $\{q_{sk}(x,t)\}$ be 
differentiable with respect to  $t$ and let coefficients $\{Q_k(x,t)$ and $\{Q_{sk}(x,t)\}$ be 
differentiable with respect to  $x$ on the domain $ {\cal I}_1 \times {\cal I}_2$.
Assume also that matrix functions $\{q_k(x,t)$, $\{q_{sk}(x,t)\}$,
$\{Q_k(x,t)$, $\{Q_{sk}(x,t)\}$, $\big\{\frac{\p}{\p t}q_k(x,t)\big\}$, and 
$\big\{\frac{\p}{\p t}q_{sk}(x,t)\big\}$
are continuous with respect to $x$ and $t$, and that zero curvature
equation
\eqref{0.1}, where $G$ and $F$
are given by \eqref{g16} and \eqref{g17},
holds. Then there is the fundamental solution of \eqref{g15}  normalized
by the condition
\begin{align} \label{ap1}
w(0,0,z)=I_{m}.
\end{align} 
\end{Cy}
\begin{proof}. Put
\begin{align} \label{ap2}
w(x,t,z)=W(x,t,z)R(t,z).
\end{align} 
By \eqref{0.4}, \eqref{0.5}, and \eqref{0.6'} we see that  $w$ given by \eqref{ap2}
satisfies \eqref{g15}. According to the second relations in \eqref{0.4} and \eqref{0.5}
equality \eqref{ap1} holds too.
\end{proof}

Further asume that $G$, $F$, and coefficients in \eqref{g16} and \eqref{g17}
satisfy conditions of Corollary \ref{Rk3.3}, and that $w$ is given by \eqref{ap2}.

When we deal with two auxiliary linear systems, 
we fix $n \in \BN$, three $n \times n$ parameter matrices,
namely $A_1$, $A_2$, and $S(0,0)$, and two $n\times m$ parameter matrices,
namely, $\Pi_1(0,0)$ and $\Pi_2(0,0)$. These matrices are chosen so
that they satisfy the matrix identity
\begin{equation}       \label{g18}
A_1S(0,0)-S(0,0)A_2=\Pi_1(0,0)\Pi_2(0,0)^*.
\end{equation}
Compare \eqref{g18} with a similar matrix identity   (\ref{g2}) for parameter matrices
$A_k$, $\Pi_k(0)$, and $S(0)$  in Appendix. Matrix functions $\Pi_1(x,t)$, $\Pi_2(x,t)$, and $S(x,t)$
are determined by the initial values $\Pi_1(0,0)$, $\Pi_2(0,0)$, and $S(0,0)$,
respectively, differential equations (\ref{g3})--(\ref{g5})  with respect to derivatives in  $x$
and similar equations with respect to derivatives in $t$. That is, 
$\Pi_1$, $\Pi_2$, and $S$ satisfy  equations:
\begin{align}       \label{ng3}
&\big(\Pi_1 \big)_x=\sum_{k=0}^r A_1^k \Pi_1 q_k+\sum_{s=1}^l \sum_{k=1}^{r_s}(A_1 -c_s I_n)^{-k}
\Pi_1 q_{sk},
\end{align}
\begin{align}
      \label{ng4}
&\big(\Pi_2^* \big)_x=-\Big(\sum_{k=0}^r  q_k\Pi_2^*A_2^k+\sum_{s=1}^l \sum_{k=1}^{r_s}q_{sk}\Pi_2^*(A_2 -c_s I_n)^{-k}
\Big),
\\   \label{ng5}
S_x=&\sum_{k=1}^r\sum_{j=1}^k A_1^{k-j} \Pi_1 q_k\Pi_2^*A_2^{j-1}-\sum_{s=1}^l \sum_{k=1}^{r_s}
\sum_{j=1}^k
(A_1 -c_s I_n)^{j-k-1} \\ \nonumber & \times
\Pi_1 q_{sk}\Pi_2^*(A_2 -c_s I_n)^{-j},
\end{align}
which coincide with (\ref{g3})--(\ref{g5}), and additional equations with respect to $t$:
\begin{align}       \label{g3'}
&\big(\Pi_1 \big)_t=\sum_{k=0}^R A_1^k \Pi_1 Q_k+\sum_{s=1}^L \sum_{k=1}^{R_s}(A_1 -C_s I_n)^{-k}
\Pi_1 Q_{sk},\\
      \label{g4'}
&\big(\Pi_2^* \big)_t=-\Big(\sum_{k=0}^R  Q_k\Pi_2^*A_2^k+\sum_{s=1}^L \sum_{k=1}^{R_s}Q_{sk}\Pi_2^*(A_2 -C_s I_n)^{-k}
\Big),
\\     \label{g5'}
S_t=&\sum_{k=1}^R \sum_{j=1}^k A_1^{k-j} \Pi_1 Q_k\Pi_2^*A_2^{j-1}-\sum_{s=1}^L \sum_{k=1}^{R_s}
\sum_{j=1}^k
(A_1 -C_s I_n)^{j-k-1} \\ \nonumber & \times
\Pi_1 Q_{sk}\Pi_2^*(A_2 -C_s I_n)^{-j}.
\end{align}
We require
\begin{equation}       \label{r30}
\{c_s\}\cap \s(A_k)=\emptyset, \quad \{C_s\}\cap \s(A_k)=\emptyset \quad (k=1,2),
\end{equation}
where $\s(A)$ is the spectrum of $A$.
Then, Theorem \ref{GBDTrd}  provides expessions for derivatives
$\big(w_A(x,t,z)\big)_x$ and $\big(w_A(x,t,z)\big)_t$:
\begin{align}& \label{r30'}
\big(w_A\big)_x=\wt G w_A-w_AG, \quad \big(w_A\big)_t=\wt F w_A-w_A F,
\end{align}
where $\wt G$ has the same structure as $G$ and is given by formulas
\eqref{g9}-\eqref{g13}. Similarly $\wt F$ has the same structure as $F$, namely,
\begin{align}&
 \label{g43}
\wt F(x,t,z)=-\sum_{k=0}^R  z^k \wt Q_k(x,t)-\sum_{s=1}^L \sum_{k=1}^{R_s}(z-C_s)^{-k}
\wt Q_{sk}(x,t),
\end{align}
where coefficients are given by  \eqref{g10}-\eqref{g13} after substitution $Q_k$, $Q_{sk}$,
$\wt Q_k$, $\wt Q_{sk}$, $R$, $R_s$, $L$, and $L_s$ instead of
$q_k$, $q_{sk}$,
$\wt q_k$, $\wt q_{sk}$, $r$, $r_s$, $l$, and $l_s$, respectively, in those formulas.
The matrix function $w_A$ in \eqref{r30'}  has the  form
\begin{equation}       \label{g27}
w_A(x, t, z)=  I_m-\Pi_2(x,t)^*S(x,t)^{-1}(A_1-z I_n)^{-1}\Pi_1(x,t)
\end{equation}
(compare with   \eqref{g7}).
By  \eqref{g15} and \eqref{r30'} we have
\begin{equation}       \label{g35}
\wt w _x=\wt G \wt w, \quad  \wt w_t=\wt F\wt w, \quad \wt w(x,t,z):=w_A(x,t,z)w(x,t,z) .
\end{equation}

The following theorem shows that $\wt G$ and $\wt F$ satisfy zero curvature
equation
\begin{align}& \label{g42}
\wt G_t- \wt F_x+[\wt G, \wt F]= 0
\end{align}
on  the domain $\cld_S$ of the points of invertibility
of $S$:
\begin{align}& \label{g40}
\cld_S={\cal I}_1 \times {\cal I}_2 \setminus \mclo_S, \quad \mclo_S=
\{(x,t)\, :   \det S(x,t)=0\}.
\end{align}
\begin{Tm} \label{WF}
Let $G$, $F$, and coefficients in \eqref{g16} and \eqref{g17}
satisfy conditions of Corollary \ref{Rk3.3}. 
Assume that $\Pi_1$, $\Pi_2$, and $S$ satisfy \eqref{g18}-\eqref{g5'}.
Then,  the coefficients
$\{\wt q_k\}$, $\{\wt q_{sk}\}$ and $\{\wt Q_k\}$, $\{\wt Q_{sk}\}$
of  $\wt G$ and $\wt F$, respectively,  are  continuous on $\cld_S$
together with derivatives $\{ \frac{\p}{\p t}\wt q_k\}$ and $\{\frac{\p}{\p t} \wt q_{sk}\}$.
Moreover, zero curvature
equation
 \eqref{g42} holds on $\cld_S$.
\end{Tm}
\begin{proof}. By  \eqref{g10}-\eqref{g13} 
and \eqref{ng3}-\eqref{g5'} the differentiability and continuity  statements 
of our theorem for the coefficients of $\wt G$ and $\wt F$ are true. 
Moreover, it follows from \eqref{0.1} that coefficients $\{ \frac{\p}{\p x} Q_k\}$ 
and $\{\frac{\p}{\p x}  Q_{sk}\}$, and hence also coefficients
$\{ \frac{\p}{\p x}\wt Q_k\}$ and $\{\frac{\p}{\p x} \wt Q_{sk}\}$, are continuous
 on $\cld_S$. 
Therefore, the matrix functions $\wt G$, $\wt G_t$,   $\wt F$, and $\wt F_x$
are continuous on the domain $\cld_S$.

Thus, taking into account \eqref{g35} we see that
$\wt w$, $\wt w_x$, $\wt w_t$, and $\wt w_{xt}$ exist and are continuous
on $\cld_S$. Hence, the conditions of the stronger formulation
of the theorem on mixed derivative, which is already used in the proof of Theorem
\ref{TmM}, are fulfilled and $\wt w_{xt}=\wt w_{tx}$ in the interior $(\cld_S)_i$ of
$\cld_S$. Using \eqref{g35} rewrite the equality $\wt w_{xt}=\wt w_{tx}$ in the form
\begin{align}& \label{g41}
(\wt G_t+\wt G \wt F)\wt w= (\wt F_x+\wt F \wt G)\wt w.
\end{align}
It follows from \eqref{g15} and \eqref{ap1} that $w$ is invertible on ${\cal I}_1 \times {\cal I}_2 $,
and it follows from \eqref{ap3} that $w_A$ is invertible on $\cld_S$. Thus
$\wt w =w_A w$ is invertible on $\cld_S$. Now,  it follows from \eqref{g41} that
\eqref{g42} holds  in $(\cld_S)_i$. By continuity, \eqref{g42} holds  on $\cld_S$.
\end{proof}
\begin{Rk}\label{F2}
Theorem \ref{WF} is basic to construct solutions and wave functions
of  integrable equations via GBDT. The corresponding normalized wave functions $\breve w$ 
(see  \cite{GKS2, GKS6, SaA2, SaA2', SaA3, SaA4, SaA6}) have the form :
\begin{align}& \label{g44}
\breve w(x,t,z)=\wt w(x,t,z)w_A(0,0,z)^{-1}=w_A(x,t,z)w(x,t,z)w_A(0,0,z)^{-1}.
\end{align}
\end{Rk}
 Compatibility of the equations \eqref{ng3}-\eqref{g5'} is a separate question.
In full generality it will be addressed elsewhere, and here we consider an important
and characteristic example of the compatibility of equations \eqref{ng3} and \eqref{g3'},
when $G$ and $F$   are polynomials.
\begin{Pn} Let $G$ and $F$ be polynomials
\begin{align} \label{g45}&
G(x,t,z)=-\sum_{k=0}^r  z^k q_k(x,t), \quad F(x,t,z)=-\sum_{s=0}^R  z^s Q_s(x,t),
 \end{align}
such that the conditions of  Corollary \ref{Rk3.3} are satisfied.
 Then, the corresponding equations \eqref{ng3} and \eqref{g3'}, which determine  $\Pi_1$, are compatible.
\end{Pn}
\begin{proof}. When $G$ and $F$ are polynomials, equations
  \eqref{ng3} and \eqref{g3'} take the form
\begin{align} \label{g46}&
\big(\Pi_1\big)_x=\sum_{k=0}^r  A_1^k\Pi_1 q_k, \quad \big(\Pi_1\big)_t=\sum_{s=0}^R  A_1^s \Pi_1Q_s.
 \end{align}
Denote the $p$-th column of $\Pi_1$ by $\big(\Pi_1\big)_p$ and introduce a block vector
with $\big(\Pi_1\big)_p$ as its blocks:
\begin{align} \label{g47}&
\ora{\Pi}_1=\left[
\begin{array}{c}
\big(\Pi_1\big)_1\\ \cdots \\ \big(\Pi_1\big)_m
\end{array}
\right] \in \BC^{mn}.
 \end{align}
Equations \eqref{g46} can be rewritten in an equivalent form in terms of  $\ora{\Pi}_1$:
\begin{align} \label{g48}&
\big(\ora{\Pi}_1\big)_x=\g\ora{\Pi}_1 , \quad 
\big(\ora{\Pi}_1\big)_t=\G \ora{\Pi}_1, \\ \label{g49}&
\g(x.t):=\sum_{k=0}^r q_k(x,t)^T \otimes A_1^k, \quad \G(x,t):=\sum_{s=0}^R Q_s(x,t)^T \otimes A_1^s,
 \end{align}
where $q^T$ denotes the transpose of a matrix $q$, and $q\otimes A$ is the Kronecker product
of matrices $q$ and $A$. As $\{q_k\}$ and $\{Q_k\}$ satisfy conditions of  Corollary \ref{Rk3.3},
so $\g$ and $\G$ satisfy the differentiability and continuity conditions of Theorem \ref{TmM},
and to prove the compatibility it remains only to show that zero curvature equation
\begin{align} \label{g50}&
\g_t-\G_x+[G,f]=0
 \end{align}
holds. For that purpose consider  the block in the $i$-th block row and in the $p$-th block column
in \eqref{g50}.  We get an equality
\begin{align} \nn&
\sum_{k=0}^r \Big(\frac{\p}{\p t}q_k\Big)_{pi}A_1^k-\sum_{s=0}^R \Big(\frac{\p}{\p x}Q_s\Big)_{pi}A_1^s
+\Big(\sum_{k=0}^r  (q_k)_i^T\otimes A_1^k \Big)\Big(\sum_{s=0}^R(Q_s^T)_p\otimes A_1^s  \Big)
\\ \label{g51} &
-
\Big(\sum_{s=0}^R (Q_s)_i^T\otimes A_1^s \Big)\Big(\sum_{k=0}^r(q_k^T)_p\otimes A_1^k  \Big)=0,
 \end{align}
 where $ (q_k)_i^T$ is the transpose of the $i$-th column of $q_k$.  Rewrite \eqref{g51} as
 \begin{equation} \nn
\sum_{k=0}^r \Big(\frac{\p}{\p t}q_k\Big)_{pi}A_1^k-\sum_{s=0}^R \Big(\frac{\p}{\p x}Q_s\Big)_{pi}A_1^s
+\sum_{k=0}^r\sum_{s=0}^R\Big((Q_s q_k)_{pi}-(q_kQ_s)_{pi}\Big)A_1^{k+s}=0.
\end{equation}
It is immediate that independently from the choice of $A_1$ the equality above
follows from the equality
\begin{align} \label{g52}&
\sum_{k=0}^r \Big(\frac{\p}{\p t}q_k\Big)_{pi}z^k-\sum_{s=0}^R \Big(\frac{\p}{\p x}Q_s\Big)_{pi}z^s
+\sum_{k=0}^r\sum_{s=0}^R\Big((Q_s q_k)_{pi}-(q_kQ_s)_{pi}\Big)z^{k+s}=0.
\end{align}
In other words equation \eqref{g50} follows from
\begin{align} \label{g53}&
\sum_{k=0}^r z^k\frac{\p}{\p t}q_k-\sum_{s=0}^R z^s\frac{\p}{\p x}Q_s
+\sum_{k=0}^r\sum_{s=0}^Rz^{k+s}\Big(Q_s q_k-q_kQ_s\Big)=0.
\end{align}
Notice that in view of  \eqref{g45}  formula \eqref{g53} is equivalent to
\eqref{0.1}.
Thus, \eqref{g53} holds, and so \eqref{g50} holds too.
\end{proof}

 %%%%%%%%%%%%%%%%%%%%%%%%%%%%%%%%%%%%%%%%%%

{\bf Acknowledgement.}
The work of A.L. Sakhnovich was supported by the Austrian Science Fund (FWF) under
Grant  no. Y330.

%%%%%%%%%%%%%%%%%%%%%%%%%%%%%%%%%%%%%%%%%%%%%%
%%%%%%%%%%%%%%%%%%%%%%%%%%%%%%%%%%%%%%%%%%%%%%%%
\appendix

\section{Appendix. GBDT for system  depending rationally on
spectral  parameter}\label{RD}
\setcounter{equation}{0}

In this appendix we consider the  GBDT version of the
B\"acklund-Darboux transformation (BDT) for a general case of first order
system  depending rationally on the spectral parameter $z$:
\begin{equation}       \label{g1}
w_x=Gw, \quad G(x,z)=-\Big(\sum_{k=0}^r  z^k q_k(x)+\sum_{s=1}^l \sum_{k=1}^{r_s}(z-c_s)^{-k}
q_{sk}(x)\Big),
\end{equation}
where $x \in {\cal I}$,   and the coefficients
$q_k(x)$ and $q_{sk}(x)$ are $m \times m$ locally integrable matrix functions. To simplify notations we assume that
${\cal I}$ is either interval $[0, \,b]$ $(0<b<\infty)$ or interval $[0, \,b)$ $(0<b\leq \infty)$.
In our presentation of GBDT  we follow Section 3 of the review \cite{SaA6}. Further references one can find
in Introduction and  \cite{SaA6}.

As GBDT is a so called iterated BDT we fix an integer $n>0$. Next, we fix five  matrices,
namely, $n \times n$ matrices $A_k$ ($k=1,2$) and $S(0)$, and $n \times m$ matrices
$\Pi_k(0)$ ($k=1,2$). It is required that these matrices form an $S$-node, that is,
the identity
\begin{equation}       \label{g2}
A_1S(0)-S(0)A_2=\Pi_1(0)\Pi_2(0)^*
\end{equation}
holds. Matrix functions $\Pi_k(x)$ are introduced via initial values $\Pi_k(0)$ and linear differential
equations:
\begin{eqnarray}       \label{g3}
&&\big(\Pi_1 \big)_x=\sum_{k=0}^r A_1^k \Pi_1 q_k+\sum_{s=1}^l \sum_{k=1}^{r_s}(A_1 -c_s I_n)^{-k}
\Pi_1 q_{sk},\\
      \label{g4}
&&\big(\Pi_2^* \big)_x=-\Big(\sum_{k=0}^r  q_k\Pi_2^*A_2^k+\sum_{s=1}^l \sum_{k=1}^{r_s}q_{sk}\Pi_2^*(A_2 -c_s I_n)^{-k}
\Big),
\end{eqnarray}
where $\{q_k\}$ and $\{q_{sk}\}$ are  coefficients from $G$.
Compare  (\ref{g1}) with (\ref{g4}) to see that $\Pi_2^*$ can be viewed
as a generalized eigenfunction of the system $u_x=Gu$.

Matrix function $S(x)$ is introduced via $\frac{d}{dx}S$ by the equality
\begin{eqnarray}       \label{g5}
S_x&=&\sum_{k=1}^r\sum_{j=1}^k A_1^{k-j} \Pi_1 q_k\Pi_2^*A_2^{j-1}-\sum_{s=1}^l \sum_{k=1}^{r_s}
\sum_{j=1}^k
(A_1 -c_s I_n)^{j-k-1} \\ \nonumber && \times
\Pi_1 q_{sk}\Pi_2^*(A_2 -c_s I_n)^{-j}.
\end{eqnarray}
Equality   (\ref{g5}) is chosen so that the identity
$\Big(A_1S-SA_2\Big)_x=\Big(\Pi_1\Pi_2^*\Big)_x$ holds.
Hence, taking into account (\ref{g2}) we have
\begin{equation}       \label{g6}
A_1S(x)-S(x)A_2=\Pi_1(x)\Pi_2(x)^*, \quad x\in {\cal I}.
\end{equation}
By Theorem \ref{GBDTrd} below, the Darboux matrix for system (\ref{g1}) has the form (\ref{0.2}) :
\begin{equation}       \label{g7}
w_A(x, z)=  I_m-\Pi_2(x)^*S(x)^{-1}(A_1-zI_n)^{-1}\Pi_1(x).
\end{equation}
In other words, $w_A$ satisfies the equation
\begin{equation}       \label{g8}
\frac{d}{d x}w_A(x, z)=\wt G(x,z)w_A(x, z)-w_A(x, z)G(x,z),
\end{equation}
where $\wt G$ has the same structure as $G$:
\begin{equation}       \label{g9}
\wt G(x,z)=-\Big(\sum_{k=0}^r  z^k \wt q_k(x)+\sum_{s=1}^l \sum_{k=1}^{r_s}(z-c_s)^{-k}
\wt q_{sk}(x)\Big).
\end{equation}
The transformed coefficients $\wt q_k$ and $\wt q_{sk}$ are given by the formulas
\begin{equation}       \label{g10}
\wt q_k=q_k-\sum_{j=k+1}^r\Big(q_jY_{j-k-1}-X_{j-k-1}q_j+\sum_{i=k+2}^jX_{j-i}q_j
Y_{i-k-2}\Big),
\end{equation}
\begin{equation}
  \label{g11}
\wt q_{sk}=q_{sk}
+\sum_{j=k}^{r_s}\Big(q_{sj}Y_{s,k-j-1}-X_{s,k-j-1}q_{sj}-\sum_{i=k}^jX_{s,i-j-1}q_{sj}
Y_{s,k-i-1}\Big),
\end{equation}
where $X_k(x)$, $Y_k(x)$, $X_{sk}(x)$, and $Y_{sk}(x)$ are expressed in terms of the
matrices $A_k$ and matrix functions $S(x)$ and $\Pi_k(x)$:
\begin{eqnarray}       \label{g12}
&&X_k=\Pi_2^*S^{-1}A_1^k\Pi_1, \quad X_{sk}=\Pi_2^*S^{-1}(A_1-c_s I_n)^k\Pi_1, \\
&&  \label{g13}
Y_k=\Pi_2^*A_2^kS^{-1}\Pi_1,  \quad Y_{sk}=\Pi_2^*(A_2-c_s I_n)^kS^{-1}\Pi_1.
\end{eqnarray}
Denote the spectrum of matrix $A$ by $\s(A)$.
\begin{Tm} \label{GBDTrd} \cite{SaA3} Let first order system   (\ref{g1}) and five
matrices $S(0)$, $A_k$, and $\Pi_k$  $(k=1,2)$ be given. Assume that
the identity  (\ref{g2}) holds and that $\{c_s\}\cap \s(A_k)=\emptyset$ $(k=1,2)$. Then, in the
points of invertibility of $S$, the transfer matrix function $w_A$ given by
(\ref{g7}), where $S$ and $\Pi_k$ are determined by (\ref{g3})--(\ref{g5}),
 satisfies equation (\ref{g8}), where $\wt G$ is determined
 by the formulas (\ref{g9})--(\ref{g13}).
 \end{Tm}
\begin{Rk} \label{det}
The matrix function $w_A$ is invertible, since it can be derived from 
\eqref{g6} that 
\begin{equation}       \label{ap3}
w_A(x, z)^{-1}=  I_m+\Pi_2(x)^*(A_2-zI_n)^{-1}S(x)^{-1}\Pi_1(x).
\end{equation}
\end{Rk}

%%%%%%%%%%%%%%%%%%%%%%%%%%%%%%%%%%%%%%%%%%%
%%%%%%%%%%%%%%%%%%%%%%%%%%%%%%%%%%%%%%%%%%%%

\end{document}